\documentclass[conference]{IEEEtran}
\IEEEoverridecommandlockouts

\usepackage{amsmath,amssymb,amsfonts}
\usepackage{booktabs}
\usepackage{cite}
\usepackage{graphicx}
\usepackage{siunitx}

\newcommand*{\rtol}{\text{rtol}}
\newcommand*{\atol}{\text{atol}}
\newcommand*{\wrms}[1]{\|#1\|_{\text{WRMS}}}

\begin{document}

\title{A Study of Fixed versus Adaptive Step Sizes for Phasor Dynamics Simulations
\thanks{%
This work was performed under the auspices of the U.S.\ Department of Energy by Lawrence Livermore National Laboratory under Contract DE-AC52-07NA27344. LLNL-CONF-2026903.
This work was supported by the U.S.\ Department of Energy, Office of Advanced Scientific Computing Research and by the Office of Electricity through the Scientific Discovery through Advanced Computing (SciDAC) program, including support from the Frameworks, Algorithms and Scalable Technologies for Mathematics (FASTMath) Institute and from the Next Generation Grid Simulations Partnership.
}
}

\author{\IEEEauthorblockN{Steven B. Roberts}
\IEEEauthorblockA{\textit{Center for Applied Scientific Computing} \\
\textit{Lawrence Livermore National Laboratory}\\
Livermore, California, USA \\
roberts115@llnl.gov}
\and
\IEEEauthorblockN{Luke Lowery}
\IEEEauthorblockA{\textit{Electrical \& Computer Engineering} \\
\textit{Texas A\&M University}\\
College Station, Texas, USA \\
lukel@tamu.edu}
\and
\IEEEauthorblockN{Slaven Peles}
\IEEEauthorblockA{
\textit{Oak Ridge National Laboratory}\\
Oak Ridge, Tennessee, USA \\
peless@ornl.gov}
\and
\IEEEauthorblockN{Adam Birchfield}
\IEEEauthorblockA{\textit{Electrical \& Computer Engineering} \\
\textit{Texas A\&M University}\\
College Station, Texas, USA \\
abirchfield@tamu.edu}
\and
\IEEEauthorblockN{Carol S. Woodward}
\IEEEauthorblockA{\textit{Center for Applied Scientific Computing} \\
\textit{Lawrence Livermore National Laboratory}\\
Livermore, California, USA \\
woodward6@llnl.gov}
}

\maketitle

\begin{abstract}
Transient phasor dynamics simulations model the response of power systems to disturbances using differential-algebraic equations (DAEs). Fixed-step integration is widely used to solve these DAEs numerically, but a step size suitable for resolving a fault transient may be unnecessarily small during the rest of a simulation. In this work, we provide a quantitative comparison of fixed and adaptive backward differentiation formulas on a suite of seven grid test cases ranging from 200 to 70,000 buses. Adaptive step sizes automatically shrink near faults, perturbations, and controller limits, and grow during less active portions of the simulation. At the global error achieved by a fixed quarter-cycle step size, adaptive results give 1.8 to 2.6 times speedups. Local error distributions, step size profiles, and work-precision plots show how accuracy, nonlinear solver, and step selection contribute to these differences. The results establish adaptive integration as the preferred approach for these simulations.
\end{abstract}

\begin{IEEEkeywords}
differential-algebraic equations, phasor dynamics, time integration.
\end{IEEEkeywords}

\section{Introduction}

Modeling and simulation of power systems are essential for grid operators and planners to evaluate grid reliability, the behavior following faults and other disturbances, and control strategies.
Simulations include a range of system components such as generators, buses, controllers, exciters, and stabilizers, each of which can introduce distinct dynamics.
As inverter-based resources are increasingly deployed into the grid, simulation challenges from their complexity and relatively fast time scales are becoming more apparent.

This paper considers phasor dynamics simulations that use differential-algebraic equations (DAEs) as the mathematical model.
For most systems of interest, these models cannot be solved analytically and instead are numerically approximated using a time integration method.
One key parameter which affects the simulation efficiency is the step size.
It can be fixed throughout a simulation or adapted to an accuracy criterion.

Fixed-step integration is widely used in practice, but the dynamic nature of phasor dynamics simulations suggests adaptive integration could provide performance benefits.
Adaptive time integration has been studied in power systems since the 1970s \cite{fuller1973vista}, with subsequent work addressing extended-term dynamic simulations \cite{sanchez1993extended,sanchez1995variable,astic1994mixed}, efficient integration of DAEs \cite{pessanha2013making,gao2017transient}, and strategies for adapting step size and order \cite{changgang2024uneven,huang2024pi,lindi2025adaptive}.

The main contributions of this paper are a quantitative assessment of fixed versus adaptive step sizes for a suite of phasor dynamics test cases and the articulation of computational tradeoffs in a state-of-the-art DAE solver.
Beyond performance comparisons, we seek to answer several questions about the simulation accuracy and solver behavior.
When switching from fixed-step to adaptive integration, what tolerances provide similar accuracy?
For a given fixed step size, how large are the local errors committed at each step?
How do time discretization errors and step sizes during a fault compare to those before and after the fault?
What is the distribution of step sizes when selected adaptively?

We use the open-source GridKit library \cite{hicss2027gridkit} as a simulation testbed to answer these questions.
For time integration in GridKit, we use IDA from the SUNDIALS library \cite{hindmarsh2005sundials,gardner2022enabling,roberts2026sundials}.
IDA, whose adaptivity origins derive from the DASSL/DASPK line of codes \cite{osti_5882821,ascher1998computer}, uses backward differentiation formula (BDF) methods with adaptive step size and order.

\section{Solver Considerations}

The cost and accuracy of a phasor dynamics simulation depend on three closely connected choices: the time discretization, the nonlinear solver at each step, and the strategy for selecting the step size.
To explain their interaction, consider the implicit DAE
\begin{equation} \label{eq:DAE}
    0 = F(t, y, y'),
    \quad
    y(t_0) = y_0,
    \quad
    y'(t_0) = y'_0,
\end{equation}
where $y(t) \in \mathbb{R}^{N}$.
The residual function $F$ contains differential equations for component dynamics and algebraic equations to impose network constraints like Kirchhoff's laws.
Fault events can change the equations in $F$ and may require a DAE solver to restart with the updated network configuration.

\subsection{BDF Discretization}

For a time step from $t_{n-1}$ to $t_n = t_{n-1} + h_n$, a BDF method
\cite[Section~5]{ascher1998computer} of order $q$ approximates the derivative by
\begin{equation} \label{eq:BDF}
    y'_n = \frac{1}{h_n} \sum_{i=0}^{q}\alpha_{n,i}y_{n-i}.
\end{equation}
The coefficients $\alpha_{n,i}$ depend on $q$ and the recent step sizes\footnote{IDA's fixed-leading-coefficient formulation uses a slight generalization of \eqref{eq:BDF} but ultimately leads to a nonlinear system of the form \eqref{eq:BDF_nonlinear_system}.}.
Substituting \eqref{eq:BDF} into \eqref{eq:DAE} defines a nonlinear system for $y_n$
\begin{equation} \label{eq:BDF_nonlinear_system}
    0 = F(t_n, y_n, h_n^{-1} \alpha_{n,0} y_n + \beta_n),
\end{equation}
with $\beta_n = h_n^{-1} \sum_{i=1}^{q}\alpha_{n,i}y_{n-i}$ determined by the step history.
Both fixed-step and adaptive integration must solve this system at each time step.

\subsection{Local Error and Step Size Adaptivity}

The local error is the difference between the result of a single BDF step and the exact solution at $t_n$ assuming that the previous history values are exact.
IDA estimates this error using the difference between a predicted state and a corrected solution to \eqref{eq:BDF_nonlinear_system} \cite[Section~2.3]{hindmarsh2005sundials}.
Denoting the estimated error vector by $e_n$, the solver tests the candidate step against user-specified absolute and relative tolerances (\atol{} and \rtol{}, respectively):
\begin{equation*}
    \wrms{e_n} = \sqrt{\frac{1}{N} \sum_{i=1}^{N} \left( \frac{e_{n,i}}{\atol + \rtol |y_{n,i}|} \right)^2} \leq 1.
\end{equation*}
If this condition fails to hold, the step is rejected and reattempted with a smaller step size.
Once step $n$ is accepted, $h_{n+1}$ is estimated to keep the next step's error within the tolerances.
IDA, for example, bases its step size update on
$h_{n+1} = 0.9 h_n / \sqrt[q+1]{2 \wrms{e_n}}$ with additional heuristics for order adaptivity and bounds on the changes.

\subsection{Nonlinear Solver}

We consider variants of Newton's method for the solution of \eqref{eq:BDF_nonlinear_system}.
Each Newton iteration requires solving a linear system with (an approximation of) the Jacobian matrix
\begin{equation*} \label{eq:Jacobian}
    J_n = F_y + h_n^{-1} \alpha_{n,0} F_{y'}.
\end{equation*}
For phasor dynamics simulations, sparse direct linear solvers are well-suited.
Once an initial matrix factorization of $J_n$ is performed, it can be used to solve $J_n \delta = r_i$ efficiently for multiple vectors $r_i$.
DAE solvers exploit this by reusing a Jacobian factorization across Newton iterations and time steps \cite[Section~2.3]{hindmarsh2005sundials}.
When convergence of the nonlinear solver stalls, say after a fault, the Jacobian is updated and factorized again.
As $J_n$ depends on $h_n$, a sufficiently large change in step size also invalidates the factorization.
To mitigate this effect, the step size is only increased if it will exceed some threshold (double the previous step size for IDA).

If Newton's method fails to converge, a standard response is reattempting the step with a smaller step size.
Of course, this is only feasible when the step size is allowed to adapt.
When using a fixed step size, we rely on Armijo backtracking line search \cite[Chapter~8]{kelley1995iterative} and a larger maximum number of iterations to make Newton's method more robust.

\subsection{Summary of Tradeoffs}

Fixed-step integration is simple to implement and provides a predictable number of steps.
The biggest computational advantage is that it removes one source of variation in the Newton matrix so that its factorization can be reused longer.
However, a step size small enough to resolve a brief fault transient can be unnecessarily restrictive during smoother portions of the simulation, whereas a larger step can increase both discretization error and nonlinear convergence difficulties.
Without step reduction as a recovery mechanism, the burden of convergence falls entirely on the nonlinear solver.
Moreover, it can be difficult to select an appropriate convergence tolerance for the nonlinear solver.
If the tolerance is too small, the nonlinear system may be ``oversolved'' without improving the temporal error.
Conversely, a tolerance too large pollutes the temporal error with nonlinear solver error.

Adaptive integration resolves transients with small steps, takes larger steps during smoother intervals, and can reduce the step size after nonlinear convergence failures.
A principled nonlinear solver tolerance can be based on the user-specified \atol{} and \rtol{}.
These advantages come with error-estimation and step-selection overhead, rejected steps that incur work without advancing the solution, and possible reductions in factorization reuse.
Finally, we note that local error control does not guarantee a corresponding global error bound.

\section{Numerical Results}

We evaluate BDF integration with fixed and adaptive step sizes on seven synthetic grid test cases listed in Table~\ref{tab:case_properties} ranging from 200 to 70,000 buses~\cite{birchfield2016grid, tamurepo2026, yuan2020developing}.
All start with a steady state solution at $t = 0\si{s}$, a fault is applied from $t = 1\si{s}$ to $t = 1.15\si{s}$, and the simulation is run to $t = 10\si{s}$.

\begin{table}[htb!]
    \centering
    \caption{Seven test cases for the numerical experiments. The number of buses in a test case appears at the end of its name.}
    \label{tab:case_properties}
    \begin{tabular}{l|lrr}
        Test case & Geographic region & Diff. vars & Alg. vars \\ \hline
        ACTIVSg200 & Central Illinois & 400 & 1,138 \\
        WECC240 & Western USA & 1,804 & 3,864 \\
        ACTIVSg500 & South Carolina & 600 & 2,204 \\
        ACTIVSg2000 & Texas & 6,056 & 14,256 \\
        ACTIVSg10k & Western USA & 27,391 & 65,443 \\
        ACTIVSg25k & Northeastern USA & 64,419 & 159,609 \\
        ACTIVSg70k & Eastern USA & 138,228 & 409,984
    \end{tabular}
\end{table}

For the IDA integrator used in the numerical experiments, we set $\atol{}=10^{-6}$ and the maximum order at two.
Higher orders are not commonly used in power grid modeling, complicate the fixed-step starting procedure, and could conflate the effects of adaptive step size and order.
We configure IDA to use KLU from SuiteSparse \cite{davis2010algorithm} as the sparse direct linear solver.
For experiments with adaptive step sizes, we use several values for \rtol{} but otherwise leave the remaining adaptivity parameters and heuristics at default values.
When using a fixed step size, the relative tolerance for the nonlinear solver is set to $10^{-4}$ unless otherwise specified.
A fixed step size of one-quarter cycle ($h_n = \frac{1}{240}\si{s}$), which is standard in PowerWorld for example, serves as a baseline for comparisons.

When we compute a global or local error, we compare a numerical solution to a reference solution generated by IDA with $\rtol{} = \atol{} = 10^{-12}$.
Local error reference solutions are generated independently at each time step using the previous numerical solution, $y_{n-1}$ and $y'_{n-1}$, as the initial condition.

\subsection{Distribution of Adaptive Step Sizes}

\begin{figure}[htb!]
    \centering
    \includegraphics[width=\linewidth]{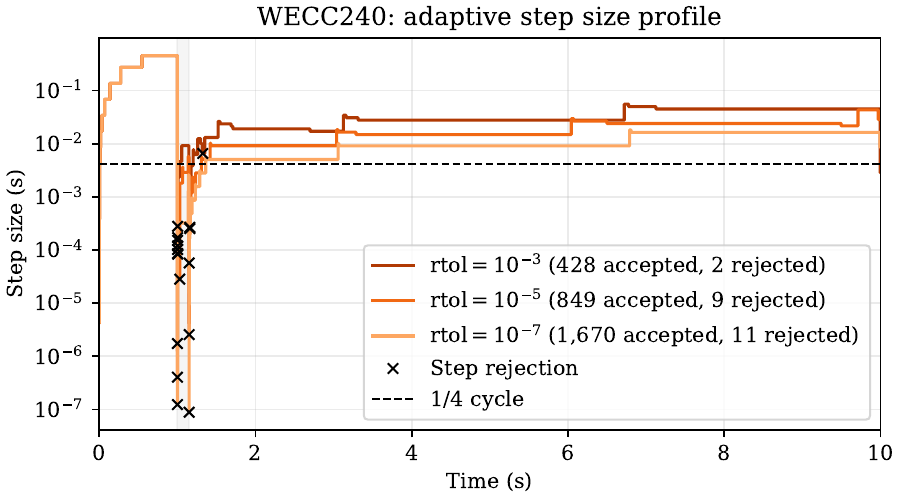}
    \caption{Adaptive step size profiles for WECC240 at three \rtol{} values.}
    \label{fig:time_step_profile}
\end{figure}

As a first demonstration of the adaptivity behavior and the effect of the \rtol{} parameter, Fig.~\ref{fig:time_step_profile} shows the evolution of the step size over the simulation time for the WECC240 test case.
Across the three values of \rtol{}, we see exponential step size growth until the fault, where it is reduced by several orders of magnitude.
The fault clearing causes a similar reduction before returning to values well above a quarter-cycle.
The number of rejected steps is insignificant at less than 4\%.
As expected, tighter tolerances generally lead to smaller step sizes.

\begin{figure}[htb!]
    \centering
    \includegraphics[width=\linewidth]{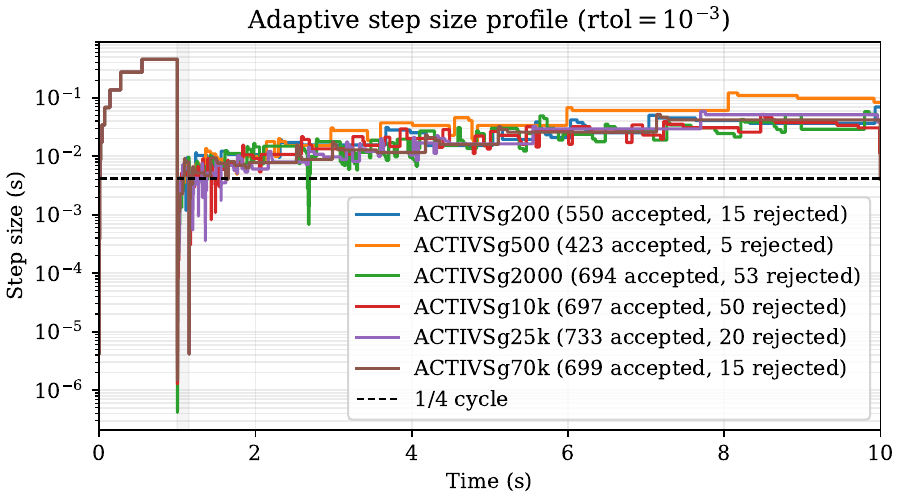}
    \caption{Adaptive step size profiles for the ACTIVSg test cases at $\rtol{} = 10^{-3}$.}
    \label{fig:aggregate_time_step_profile}
\end{figure}

\begin{figure}[htb!]
    \centering
    \includegraphics[width=\linewidth]{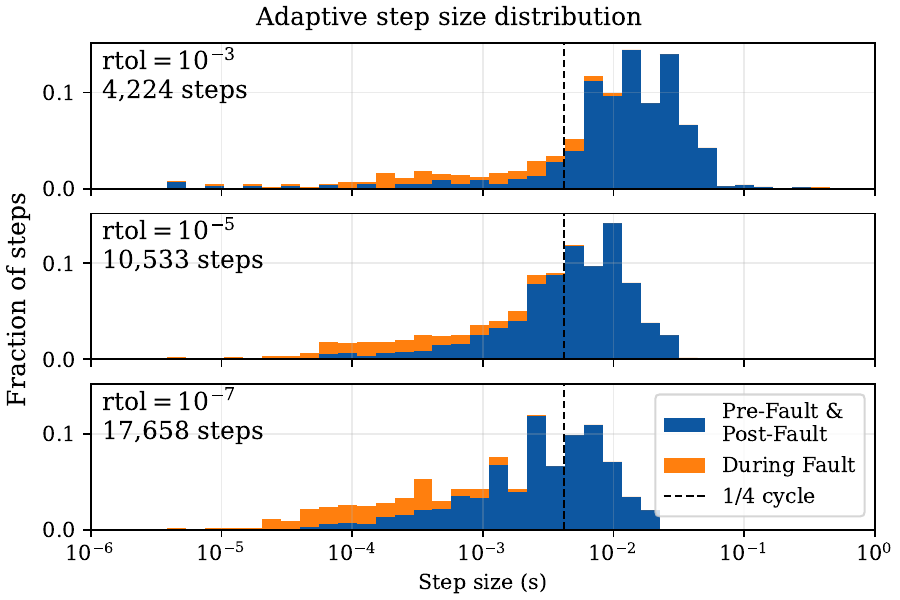}
    \caption{Adaptive step size distribution across all seven test cases at three $\rtol{}$ values. Colors distinguish steps ending during and outside the fault interval.}
    \label{fig:aggregate_time_step_histogram}
\end{figure}

The ACTIVSg test cases show similar step size reductions around the fault followed by growth during smoother intervals as shown in Fig.~\ref{fig:aggregate_time_step_profile} for $\rtol{} = 10^{-3}$.
Despite their substantially different state dimensions, the test cases exhibit similar step size profiles.
Fig.~\ref{fig:aggregate_time_step_histogram} summarizes the step sizes across all seven test cases.
For $\rtol{} = 10^{-3}, 10^{-5}, 10^{-7}$ the median step sizes are approximately $0.67$, $0.26$, and $0.14$ cycles, respectively.

\subsection{Distribution of Local Errors}

We next examine how the error introduced at each time step varies with the system dynamics.
In the visualizations, error is measured in the relative 2-norm $\|y_n-y_n^{\mathrm{ref}}\|_2/\|y_n^{\mathrm{ref}}\|_2$, where the reference solution is advanced from the preceding numerical state as described above.
This procedure measures the discrepancy over one step rather than the error accumulated since the start of the simulation.

\begin{figure}[htb!]
    \centering
    \includegraphics[width=\linewidth]{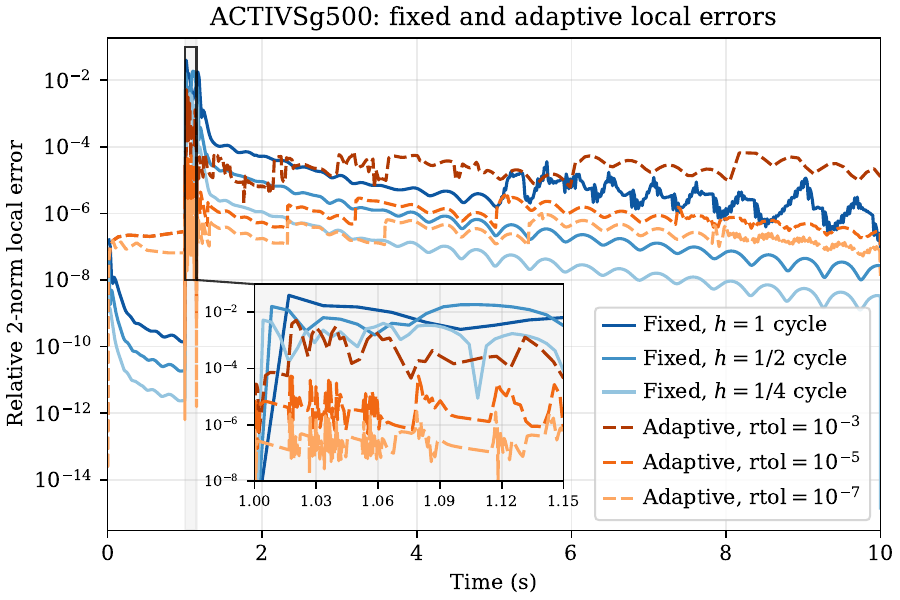}
    \caption{Relative 2-norm local errors for ACTIVSg500 using fixed-step integration with three step sizes and adaptive integration with three \rtol{} values. The inset enlarges the fault interval.}
    \label{fig:local_error}
\end{figure}

Before the fault, the fixed-step local errors are very small because the system remains close to its initial steady state.
Fault application increases these errors by several orders of magnitude, and the subsequent transient response produces substantial variation even after the fault has cleared.
Fig.~\ref{fig:local_error} compares local errors for the ACTIVSg500 test case.
The inset shows that tighter adaptive tolerances generally reduce the local errors during the fault, whereas a fixed step size cannot respond to the abrupt change in dynamics.

\begin{figure}[htb!]
    \centering
    \includegraphics[width=\linewidth]{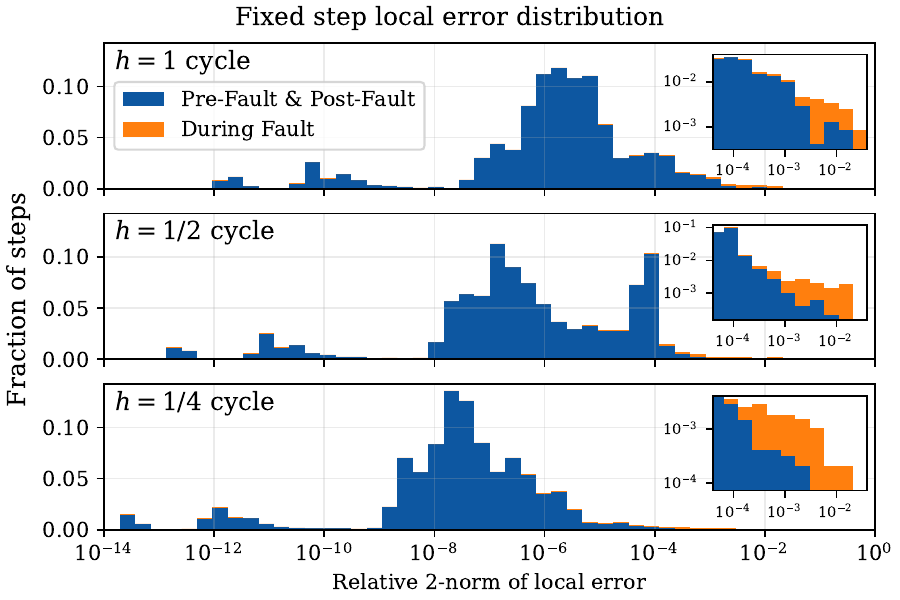}
    \caption{Aggregated relative 2-norm local error distributions for three fixed step sizes. Insets show the large-error tails on a logarithmic vertical scale.}
    \label{fig:aggregate_fixed_local_error_histogram}
\end{figure}

Fig.~\ref{fig:aggregate_fixed_local_error_histogram} shows the fixed-step local error distributions for the first four test cases in Table~\ref{tab:case_properties}.
They span over ten orders of magnitude, reflecting the difference in steady-state, fault, and post-fault dynamics.
Reducing the fixed step size shifts the bulk of the distribution toward smaller errors, but even the quarter-cycle step retains a tail around 1\% relative error.
The insets show the disproportionate contribution of steps during the fault to the tail.

\begin{figure}[htb!]
    \centering
    \includegraphics[width=\linewidth]{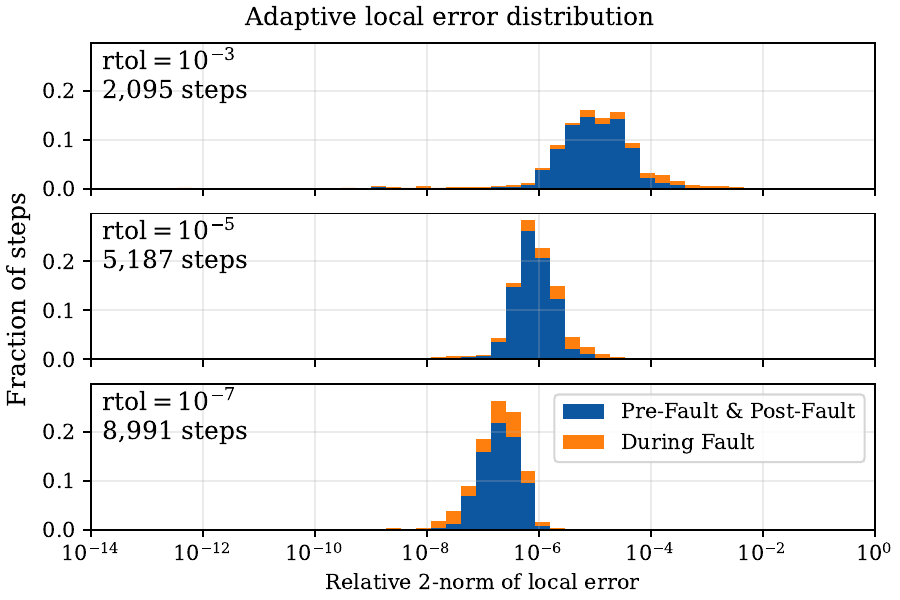}
    \caption{Aggregated relative 2-norm local error distributions for adaptive integration with three relative tolerances.}
    \label{fig:aggregate_adaptive_local_error_histogram}
\end{figure}

The corresponding adaptive distributions in Fig.~\ref{fig:aggregate_adaptive_local_error_histogram} are more concentrated, and tightening $\rtol{}$ shifts them toward smaller errors.
Fault and non-fault samples overlap more closely than in the fixed-step distributions, consistent with step size reduction compensating for the faster fault dynamics.
Together, the distributions demonstrate the difficulty of selecting one fixed step size that avoids both unnecessarily small errors during smooth intervals and large errors during disturbances.

\subsection{Performance Results}

We compare the computational cost of fixed-step and adaptive integration at comparable global accuracy, since equal tolerances do not ensure equal trajectory errors.
For these comparisons, we collect the recorded output vectors into trajectory matrices $Z=[\mathbf{z}_0,\ldots,\mathbf{z}_M]$ and $Z^{\mathrm{ref}}=[\mathbf{z}_0^{\mathrm{ref}},\ldots,\mathbf{z}_M^{\mathrm{ref}}]$ at quarter-cycle output times and define
\begin{equation}
    E_{\mathrm{traj}} = \|Z-Z^{\mathrm{ref}}\|_F / \|Z^{\mathrm{ref}}\|_F.
    \label{eq:trajectory_error}
\end{equation}
Here, $\|\cdot\|_F$ is the Frobenius norm.
Runtimes are taken from corresponding runs with monitoring and logging disabled to exclude the cost of output and reference diagnostics.

Table~\ref{tab:speedup} compares the six ACTIVSg test cases at the global error achieved by the fixed quarter-cycle baseline.
For each test case, adaptive work and relative tolerance are estimated by log-log interpolation between the two nearest adaptive error measurements that bracket the baseline error.
Speedup is the fixed work divided by the corresponding adaptive estimate.

At the baseline accuracy, adaptive integration gives runtime speedups of 1.76 to 2.64 times, reducing residual evaluations by a factor of 2.51 to 4.48 despite requiring 1.51 to 2.63 times as many matrix factorizations.
These gains are consistent with larger steps during smoother dynamics offsetting the costs of adaptivity and reduced factorization reuse.
The corresponding adaptive relative tolerances, $2.25\times10^{-5}$ to $2.51\times10^{-4}$, provide a rough starting point for matching quarter-cycle accuracy, but should be validated for each model and disturbance.

\begin{table*}[htb!]
    \centering
    \caption{Comparison of solver statistics representing the work of fixed-step and adaptive integration.
    Values in the Speedup column greater than one favor adaptive integration.
    Fixed-step integration with a quarter-cycle step size and adaptive integration based on the relative tolerance in the third column achieve approximately the global error in the second column.
    }
    \label{tab:speedup}
    \begin{tabular}{l*{11}{r}}
        \toprule
        & & & \multicolumn{3}{c}{Runtime}
        & \multicolumn{3}{c}{Residual evaluations}
        & \multicolumn{3}{c}{Matrix factorizations} \\
        \cmidrule(lr){4-6}
        \cmidrule(lr){7-9}
        \cmidrule(lr){10-12}
        Test case & Global error & Adaptive $\rtol$ & Fixed & Adaptive & Speedup & Fixed & Adaptive & Speedup & Fixed & Adaptive & Speedup \\
        \midrule
        ACTIVSg200 & $1.19\times10^{-2}$ & $1.63\times10^{-4}$ & 0.224\si{s} & 0.0982\si{s} & 2.28 & 4,907 & 1,212 & 4.05 & 30 & 95 & 0.32 \\
        ACTIVSg500 & $5.24\times10^{-3}$ & $2.82\times10^{-5}$ & 0.55\si{s} & 0.209\si{s} & 2.64 & 4,902 & 1,094 & 4.48 & 35 & 87 & 0.40 \\
        ACTIVSg2000 & $1.09\times10^{-2}$ & $1.99\times10^{-4}$ & 6.85\si{s} & 3.35\si{s} & 2.04 & 5,106 & 1,347 & 3.79 & 92 & 182 & 0.51 \\
        ACTIVSg10k & $9.44\times10^{-3}$ & $2.25\times10^{-5}$ & 55.4\si{s} & 31.5\si{s} & 1.76 & 5,086 & 2,023 & 2.51 & 81 & 213 & 0.38 \\
        ACTIVSg25k & $1.89\times10^{-2}$ & $2.51\times10^{-4}$ & 132\si{s} & 56.7\si{s} & 2.33 & 5,016 & 1,394 & 3.60 & 70 & 171 & 0.41 \\
        ACTIVSg70k & $3.54\times10^{-2}$ & $5.38\times10^{-5}$ & 347\si{s} & 155\si{s} & 2.24 & 5,145 & 1,549 & 3.32 & 106 & 160 & 0.66 \\
        \bottomrule
    \end{tabular}
\end{table*}

\begin{figure}[htb!]
    \centering
    \includegraphics[width=\linewidth]{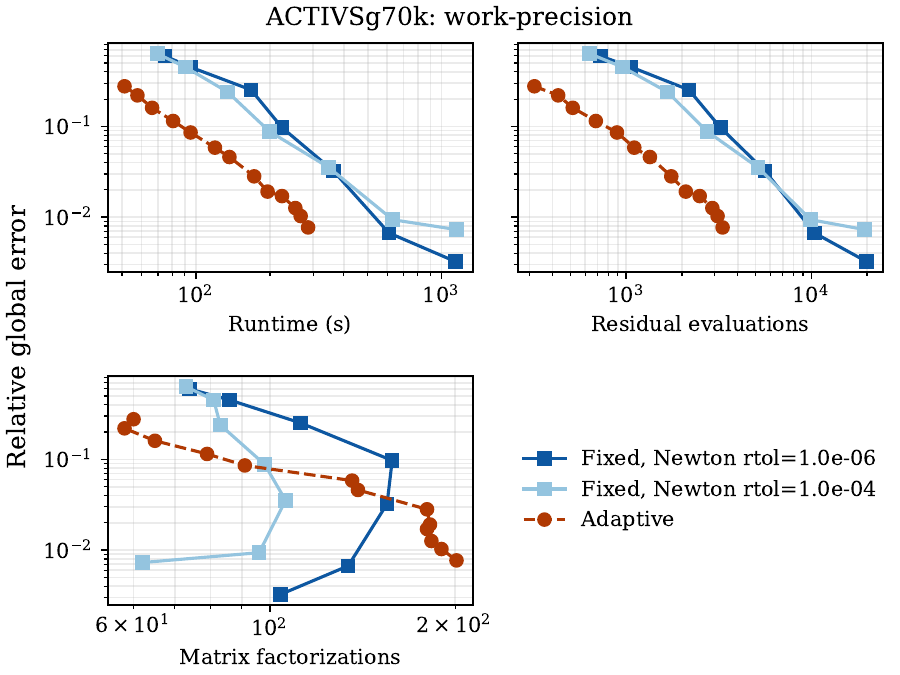}
    \caption{ACTIVSg70k work-precision curves using trajectory error~\eqref{eq:trajectory_error} and runtime (top left), residual evaluations (top right), and matrix factorizations (bottom left). Fixed-step curves vary $h_n$ at two nonlinear relative tolerances; the adaptive curve varies $\rtol{}$.}
    \label{fig:work_precision}
\end{figure}

Fig.~\ref{fig:work_precision} extends the comparison in Table~\ref{tab:speedup} across a range of accuracies for ACTIVSg70k, the largest test case.
Fixed-step runs use $h_n = \frac{2^{i}}{60}\si{s}$ for $i = -4, -3, \dots, 2$ at nonlinear solver relative tolerances of $10^{-4}$ and $10^{-6}$.
Adaptive runs use $\rtol{}=10^{-i/2}$ for $i=2,3,\ldots,14$.
Adaptive integration consistently achieves lower runtime and fewer residual evaluations.
At the finest fixed steps, tightening the nonlinear tolerance improves accuracy, indicating that nonlinear solver error can limit the benefit of reducing the step size.

\section{Conclusion}

In this work, we presented a comparison of fixed-step and adaptive integration of DAEs arising in phasor dynamics simulations.
Across the six ACTIVSg test cases, adaptive integration reduced simulation cost at comparable trajectory accuracy, giving 1.8 to 2.6 times speedups relative to a fixed quarter-cycle step size.
This is despite additional matrix factorizations and rejected attempts, indicating that the savings from adapting to changing dynamics outweigh these costs.
Our local error diagnostics further demonstrate the limitations of fixed-step integration: a step size that is adequate during smooth dynamics can introduce substantially larger errors during a disturbance.

We recommend adaptive integration as the default choice for phasor dynamics simulations.
Although this study uses second order BDF methods, we expect the benefits of adaptive step sizes to extend to other implicit methods such as the implicit trapezoidal method.
However, more comparisons are needed to assess that generality.


\bibliographystyle{IEEEtran}
\bibliography{main}

\end{document}